\documentclass[12pt]{article}

\usepackage{latexsym}
\usepackage{amssymb}
\usepackage{amsmath}
\usepackage{amsfonts}
\usepackage{mathrsfs}
\usepackage[all]{xy}
\usepackage{texdraw}
\usepackage{graphicx}
\usepackage{psfrag}
\usepackage{makeidx}

\newtheorem{Theory}{Theory}[section] 
\newtheorem{theorem}[Theory]{Theorem}
\newtheorem{lemma}[Theory]{Lemma}
\newtheorem{technicalLemma}[Theory]{Technical Lemma}
\newtheorem{corollary}[Theory]{Corollary}
\newtheorem{proposition}[Theory]{Proposition}
\newtheorem{fact}{Fact}  
\newtheorem{remark}[Theory]{Remark} 
\newtheorem{question}{Question} 
\newtheorem{conjecture}[question]{Conjecture}

\newtheorem{Ntn}{Description} 
\newtheorem{Dn}[Ntn]{Definition}

\newcommand{\be}{\begin{enumerate}}
\newcommand{\ee}{\end{enumerate}}
\newcommand{\bq}{\begin{question}}
\newcommand{\eq}{\end{question}}
\newcommand{\bcj}{\begin{conjecture}}
\newcommand{\ecj}{\end{conjecture}}
\newcommand{\bc}{\begin{corollary}}
\newcommand{\ec}{\end{corollary}}
\newcommand{\bl}{\begin{lemma}}
\newcommand{\el}{\end{lemma}}
\newcommand{\btl}{\begin{technicalLemma}}
\newcommand{\etl}{\end{technicalLemma}}
\newcommand{\bt}{\begin{theorem}}
\newcommand{\et}{\end{theorem}}
\newcommand{\bp}{\begin{proposition}}
\newcommand{\ep}{\end{proposition}}
\newcommand{\bft}{\begin{fact}}
\newcommand{\eft}{\end{fact}}
\newcommand{\brk}{\begin{remark}}
\newcommand{\erk}{\end{remark}}
\newcommand{\bd}{\begin{Dn}}
\newcommand{\ed}{\end{Dn}}

\newcommand{\ga}{\alpha}
\newcommand{\gb}{\beta}

\newcommand{\ploi}{PL_o(I) }

\newcommand{\Z}{ \mathbf Z }
 
\newcommand{\N}{\mathbf N }
 
\newcommand{\R}{\mathbf R }

\newcommand{\bSeq}[4]{\left\{#1_{#2}\right\}_{#2=#3}^{#4}}

\newcommand{\fix}{\operatorname{Fix}}
\newcommand{\supp}{\operatorname{Supp}}

\author{Collin Bleak}

\usepackage{epsfig}

\begin{document}
\centerline{\LARGE A Geometric Classification of Some Solvable Groups}
\centerline{\LARGE of Homeomorphisms}
\vspace{.3 in}

\centerline{\textbf{Abstract}} We investigate subgroups of the group
$\ploi$ of piecewise-linear, orientation-preserving homeomorphisms of
the unit interval with finitely many breaks in slope, under the
operation of composition, and also subgroups of generalized Thompson
groups $F_n$.  We find geometric criteria determining the derived
length of any such group, and use this criteria to produce a geometric
classification of the solvable and non-solvable subgroups of $\ploi$
and of the $F_n$.  We also show that any standard restricted wreath product
$C\wr T$ (of non-trivial groups) that embeds in $\ploi$ or $F_n$ must
have $T\cong \Z$.

\tableofcontents

\newpage

\setcounter{page}{1}
\section{Introduction}
  Let $\ploi$ represent the group of orientation-preserving
piecewise-linear homeomorpisms of the unit interval admitting finitely
many breaks in slope under the operation of composition.  We provide a
``geometric'' classification of the solvable subgroups of $\ploi$.  We
also show that any standard restricted wreath product $C\wr T$ (of
non-trivial groups) that embeds in $\ploi$ must have $T\cong \Z$.
This was previously unknown, and answers a question suggested by
Sapir.  All of these results apply equally well when we replace the
group $\ploi$ by any of the generalized Thompson groups $F_n$, which
can all be realized as particular subgroups of $\ploi$.  (The groups
$F_n$ were introduced by Brown in \cite{BrownFinite}, where they were
denoted $F_{n,\infty}$. These groups were later extensively studied by
Stein in \cite{SteinPLGroups}, by Brin and Guzm\'an in
\cite{BGAutomorphisms} and by Burillo, Cleary, and Stein in
\cite{BCS}.)

We proceed by creating some new objects and techniques associated with
subgroups of $\ploi$.  The new objects have a visual context, and they
help to understand the dynamics of the action of $\ploi$ on the unit
interval.  One of these objects is a special set called a
\emph{tower}, and subgroups of $\ploi$ may have many towers associated
with them.  If we call the supremum of the cardinalities of the towers
associated with a subgroup $G$ of $\ploi$ the \emph{depth of $G$},
then we can state our main result.

\bt 
\label{geoClassification}

Let $G$ be a subgroup of $\ploi$, and let $n$ be a natural number.
$G$ has derived length $n$ if and only if the depth of $G$ is $n$.
\et

Essentially, we introduce a geometric invariant of subgroups of
$\ploi$ which detects exactly the derived length of any such subgroup.
Detailed definitions are given in the next section.

Theorem \ref{geoClassification} is central to the proofs of the theory
in \cite{BPASC} providing two distinct algebraic classifications of
the solvable subgroups of $\ploi$.  The following is the more concrete
classification proven there.  First, let the symbol $A\wr B$ represent
the standard restricted wreath product as discussed in Peter Neumann's
paper \cite{NeumannW}, or in the book \cite{Meldrum} by
J. D. P. Meldrum.  Now, define $G_0=1$, the trivial group, and
inductively define the group $G_i$ by the rule $G_i =
\bigoplus_{j\in\Z}(G_{i-1}\wr\Z)$ for each positive integer $i$.  With
these definitions in place, we can state a version of one of the main
results of \cite{BPASC}; $H$ is a solvable subgroup of $\ploi$ with
derived length at most $i$ if and only if $H$ is isomorphic to a
subgroup of $G_i$.

Also, Theorem \ref{geoClassification} is central to the proofs
of the theory in \cite{BPANSC} providing an algebraic structure
theorem describing the non-solvable subgroups of $\ploi$.  Define $W_0
= 1$, the trivial group.  Now define $W_i = W_{i-1}\wr\Z$ for all
positive integers $i$.  Finally, define $W = \bigoplus_{i\in\N}W_i$.
The paper \cite{BPANSC} shows that a subgroup $H$ in $\ploi$ is
non-solvable if and only if $H$ admits a subgroup which is isomorphic
to $W$.

The idea of applying geometric analysis to sets of elements of $\ploi$ in
order to derive information about the structure of subgroups of $\ploi$ is not
new.  Brin in \cite{BrinU} finds a geometric condition under which a
subgroup of $\ploi$ is guaranteed to contain a subgroup isomorphic to
R. Thompson's group $F$ (We shall refer to this geometric condition as
Brin's ubiquity condition.  Note also that $F$ here refers to the
generalized Thompson group $F_2$ mentioned above.)  Guba (unpublished,
\cite{GubaCU}) has an argument providing a converse statement to
Brin's main result in \cite{BrinU}.

The works just mentioned set the stage for the work here, analyzing
the solvable subgroups of $\ploi$.  In particular, Thompson's group
$F$ is not in $EG$, the elementary amenable groups (a proof of this is
in Cannon, Floyd, and Parry's introduction \cite{CFP} to Thompson's
group $F$), so any subgroup of $\ploi$ that contains a copy of $F$
cannot be solvable.  This implies that in order to classify the
solvable subgroups of $\ploi$, one should consider the groups which do
not admit Brin's ubiquity condition.  It turns out that to be
solvable, subgroups of $\ploi$ have to satisfy a stricter geometric
condition than not meeting Brin's ubiquity condition (this was
known already, for example, Brin in \cite{BrinEG} finds many
non-solvable elementary amenable groups).  Describing this stronger
condition is the main work of this paper.

Our second result focusses on the types of standard restricted
wreath products (of groups) that can be found embedded in $\ploi$.  We
note that the importance of wreath products in the analysis of
subgroups of $\ploi$ has been known for a long time.  Holland in
\cite{HollandLOG} recognizes the natural occurence of many wreath
products with $\Z$ in subgroups of $Homeo_o(I)$.  Guba and Sapir in
\cite{GSFSubgroups} show that non-commutative subgroups of $F$ must
contain subgroups isomorphic to $\Z\wr\Z$.  Brin in \cite{BrinEG}
devotes a large portion of the paper to recognizing wreath products in
$\ploi$.

Our second chief result is stated below.

\bt 
\label{wreathTopClassification}

If $H$ is a subgroup of $\ploi$, and $H$ is isomorphic to a
standard restricted wreath product $C\wr T$ of nontrivial groups, then
$T\cong \Z$.  

\et

Theorem \ref{wreathTopClassification} is an easy consequence of the
lemma that $\Z\wr(\Z\times\Z)$ does not embed in $\ploi$.  While this
lemma and the resulting theorem are new, the argument we give in this
paper is not the only one which can produce the lemma.  It is an easy
fact that for each $n\in\N$, $G_n$ embeds as a subgroup of $W_{3n-1}$.
By the main result of \cite{BPASC}, if $\Z\wr(\Z\times\Z)$ embeds in
$\ploi$, then it embeds in $G_2$, and therefore in $W_5$.
A. Yu. Ol'shanskii has an algebraic proof that $\Z\wr(\Z\times\Z)$
does not embed in $W_i$ for any $i\in\N$.  Therefore, Ol'shanskii's
proof and the main result of \cite{BPASC} are also enough to show that
$\Z\wr(\Z\times\Z)$ does not embed in $\ploi$.

The author would like to thank Matt Brin for asking that he interpret
the main geometric result of this paper into algebraic language, which
has led to the results in \cite{BPASC}.  He would also like to thank
Victor Guba, A. Yu. Ol'shanskii, and Mark Sapir for various
conversations about and leading to the second main result of this
paper.

Portions of this paper contain portions of the author's dissertation
written at Binghamton University.

\subsection{Geometric concepts}
The standard definitions referring to the geometry of elements of
$\ploi$ will be reserved for a later section.  Here we will only give
the definitions required to understand the statements of our
key results.

For the rest of this section, fix $G\leq\ploi$, a subgroup of $\ploi$,
for reference.

We will call the convex hull of an orbit of a point in $I$ under the
action of $G$ an \emph{orbital of $G$}, if this set has cardinality
greater than one.  We note that orbitals are open intervals.  If
$g\in G$, then we will refer to an orbital of the group $\langle
g\rangle$ as an orbital of $g$.  Given an interval $A = (a,b)\subset
I$, if $A$ is an orbital of $g$ for some $g\in G$, then $A$ is also an
orbital of $g^n$ for any non-zero integer $n$.  To be more specific,
we sometimes associate the interval $A$ with a specific element of $G$
which has $A$ as an orbital.  Hence we commonly form \emph{signed
orbitals of $G$}, which are pairs $(A,g)$ so that $g\in G$ and $A$ is
an orbital of $g$.  Given a signed orbital $(A,g)$ of $G$, we call $A$
the orbital of $(A,g)$, and we call $g$ the signature of $(A,g)$.

Given a set $Y$ of signed orbitals of $G$, the symbol $S_Y$ will refer
to the set of signatures of the signed orbitals in $Y$.  Similarly, the
symbol $O_Y$ will refer to the set of orbitals of the signed orbitals
of $Y$.  We note that the set of signed orbitals of $\ploi$ is a
partially ordered set under the lexicographic order induced from the
partial order on subsets of $I$ (induced by inclusion) in the
first coordinate, and the left total order of the elements of $\ploi$
in the second coordinate. (Given two distinct elements of $\ploi$,
there is a maximal $x\in [0,1)$ so that the two elements agree on the set
$[0,x]$, but one of the elements will have a larger value in its right
hand derivative at $x$.  We will call this element the ``larger''
element, inducing a total order on the elements of $\ploi$, called the
left total order of $\ploi$.)

A \emph{tower $T$ of $G$} is a set of signed orbitals which satisfies
the following two criteria.
\begin{enumerate}
\item $T$ is a chain in the partial order on the signed orbitals of $G$.
\item For any $A\in O_T$, $T$ has exactly one element of the
form $(A,g)$.
\end{enumerate}

Given a tower $T$ of $G$, if $(A,g)$, $(B,h)\in T$ then one of
$A\subset B$ and $B\subset A$ holds, with equality occuring only if
$g=h$ as well. Therefore, one can visualize a tower as a stack of
nested levels that are always getting wider as one goes ``up'' the
stack.  We define the \emph{height} of a tower as its cardinality.

While we say that towers have height, we will say that their
associated groups have depth.  Thus, we will say that the \emph{depth of
$G$} is the supremum of the cardinalities of the towers of the group.

If $G$ admits two signed orbitals $(A,g)$ and $(B,h)$ so that $A =
(a_1,a_2)$ and $B = (b_1,b_2)$, with $a_1<b_1<a_2<b_2$ then we will
say that $G$ admits a \emph{transition chain of length two}.  This
definition has a natural generalization that encompasses longer
transition chains, but we will have no need for that generality here.

If $A = (a,b)$ is an orbital of the group $G$, and $G$ has an element
$g$ which has an orbital $B = (c,d)$ so that either $c = a$ or $d =
b$, then we say that $g$ has an orbital that \emph{shares an end} with
$A$.

\subsection{Some further results} 
Here we will mention a few of the consequences of Theorem
\ref{geoClassification}.  We begin with a simple corollary that is
really just a restatement of the main theorem, together with an
application of some of the simpler properties of towers.

\bc 
A subgroup $H$ of $\ploi$ is non-solvable if and only if $H$
admits towers of arbitrary height.
\ec 

The following lemmas show that the sets of support of distinct
elements in a solvable subgroup of $\ploi$ must be very simply
related.  Because of these lemmas, we can examine a few elements of a
subgroup of $\ploi$, and very often, we can detect that the subgroup
will not be solvable.

\bl 
\label{tChainTowers}
If $H$ is a subgroup of $\ploi$, and $H$ admits transition chains
of length two, then $H$ admits infinite towers.
\el

\bl 

\label{imbalancedTChains}
Suppose $H$ is a subgroup of $\ploi$ with orbital $A$.  If $H$ has an
element $h$ which has one orbital sharing exactly one end of $A$,
while the other orbitals of $h$ do not share an end with $A$, then $H$
admits transition chains of length two.
\el

\section{The geometry of $\ploi$}
Here we will describe objects and techniques associated with subgroups
of $\ploi$.  First we will describe some of what was known, the
dynamics associated with conjugation and orbitals, and then we will
introduce some new techniques involving towers and transition chains
as defined above.
\subsection{Known geometric objects of $\ploi$}
We will use notation similar to that in Brin's paper
$\cite{BrinU}$ on the ubiquitous nature of Thompson's group $F$ in
subgroups of $\ploi$.  

Most of the statements in this section are either basic (whose proofs
are left to the reader), or are standard facts (proofs can be found in
the literature (\cite{BSPLR,picric, BrinU})).  The last two statements
will be proven here.

Orbitals as defined above can also be thought of as the components of
the support of the action of a subgroup of $\ploi$ on the unit
interval.

Given an open interval $A =(a,b)\subset\R$, where $a<b$, we will refer
to $a$ as the \emph{leading end of $A$}\index{end!leading}, and to $b$
as the \emph{trailing end of $A$}\index{end!trailing}.  If the
interval is an orbital of some group $H\in\ploi$, we will refer to
the \emph{ends of the orbital}\index{orbital!end} in the same fashion.

If $h\in{}H$ and $x\in{}\supp{h}$, we will say that \emph{$h$ moves
$x$ to the left}\index{movement!left} if $xh<x$, and we will say that
\emph{$h$ moves $x$ to the right}\index{movement!right} if $xh>x$.
Given any $x\in (0,1)$, we will say that $x$ is a
\emph{breakpoint for $h$}\index{breakpoint} if the left and right
derivatives of $h$ exist at $x$, but are not equal.  We recall that by
definition, $h$ will admit only finitely many breakpoints.  If
$\mathscr{B}_h$ represents the set of breakpoints of the element $h$,
then $(0,1)\backslash\mathscr{B}_h$ is a finite collection of open
intervals, which we will call \emph{affine components of
$h$}\index{component!affine}, and which admits a natural ``left to right''
ordering.  We shall refer to the ``first'' affine component
of $h$, or the ``second'' affine component of $h$, etc.  We sometimes
will refer to the first affine component of $h$ as the \emph{leading
affine component of $h$}\index{component!leading affine}, and to the
last affine component of the domain of $h$ as the \emph{trailing
affine component of $h$}\index{component!trailing affine}.

The following are some useful remarks.

\brk
\label{finiteOrbitals}
\label{transitiveElementOrbital}
\be

\item If $A$ is an orbital for $h\in H$, then either $xh>x$ for all
points $x$ in $A$, or $xh<x$ for all points $x$ in $A$.

\item Any element $h\in\ploi$ has only finitely many
orbitals.

\item If $h\in\ploi$ and $A=(a,b)$ is an orbital of $h$, then given any
$\epsilon>0$ and $x$ in $A$, there is an integer $n$ so that
$|xh^{-n}-a|<\epsilon$ and $|xh^n-b|<\epsilon$.
\ee
\erk

pf: 

All three points can be taken as exercises for the
reader.  The third point is a modification of Lemma 3.4 in \cite{BSPLR}.
\qquad$\diamond$

Relating to the first point above, we note that if $g\in\ploi$ and $A$
is an orbital of $g$, so that $xg<x$ for some (and therefore all)
$x\in A$, then $g^{-1}$ has $xg^{-1}>x$ for all $x\in A$.  We will
often use this property to intelligently choose a signature for a
signed orbital so that the signature moves points right throughout its
orbital.

Given an orbital $A$ of $H$ we say that \emph{$h$ realizes an end of
$A$} if some orbital of $h$ lies entirely in $A$ and shares an
end with $A$.  Note that Brin uses the word ``Approaches'' for this
concept in \cite{BrinU}, but we will use ``Approaches'' to also
indicate the direction in which $h$ moves points, as follows: we will
say that \emph{$h$ approaches the end $a$ of $A$ in $A$} if $h$ has an
orbital $B$ where $B\subset A$ and $B$ has end $a$, and $h$
moves points in $B$ towards $a$.  If $A$ is
an orbital for $H$ then we say that \emph{$h\in H$ realizes $A$} if
$A$ is also an orbital for $h$.

If $g$ and $h$ are elements of $\ploi$ and there is an interval $A =
(a,b)\subset I$ so that both $g$ and $h$ have $A$ as an orbital, then
we will say that $g$ and $h$ \emph{share the orbital $A$}. 

We will now start to analyze the effects of conjugation.  A fuller
discussion is contained in \cite{picric}.

The following is standard and is left to the reader.

\bl 

Suppose $n\in\N$ and $g$, $h\in \ploi$.  Further suppose $g$ has $n$
orbitals and $\bSeq{A}{i}{1}{n}$ represents these orbitals arranged in
left to right order. The collection of open intervals
$\bSeq{A^*}{i}{1}{n}$ defined by the rule $A_i^* = A_ih$ represents the
collection of orbitals of $g^h = h^{-1}gh$.

\el

In the setting of the above lemma, we will say that the $A_ih$  are the
\emph{induced orbitals}\index{orbital!induced} of $g^h$ from $g$ by the
action of $h$.  We might also say that the orbitals of $g^h$ are induced
from the orbitals of $g$ by the action of $h$.

The following two unrelated points are also worth pointing out:
\brk 
\label{breakpoints}
\be
\item Suppose $g$, $h\in\ploi$ and $f = gh$.  If $b$ is a breakpoint of $f$
then $b$ is a breakpoint of $g$ or $bg$ is a breakpoint of $h$.
\item Let $\ga$, $\gb\in\ploi$, then $\ga^{\gb}$ has the same leading
and trailing slopes on its orbitals as\, $\ga$ has on each of its
corresponding orbitals.
\ee
\erk

pf: 

For the first point, suppose $b$ is a breakpoint of $f$, but $b$ is not a
breakpoint of $g$.  The lefthand slope of $f$
at $b$ is the product of the lefthand slopes of $g$ at $b$
and $h$ at $bg$. Also, the righthand slope of $f$
at $b$ is the product of the righthand slopes of $g$ at $b$
and $h$ at $bg$.  Since the lefthand and righthand slopes
of $g$ at $b$ are the same, the lefthand and righthand slopes of $h$
at $bg$ must be different.  In particular, $bg$ is a breakpoint for $h$.

For the second point, let $A = (a,b)$ be an orbital of $\ga^{\gb}$,
corresponding to an orbital $A' = (a',b')$ of $\ga$.  There is
$c\in(a,b)$ so that $(a,c)$ is contained in an affine component of
$\gb^{-1}$, so that $(a',c') = (a,c)\gb^{-1}$ is in an affine
component of $\ga$, and so that $(a',c'') = (a',c')\ga$ is contained
in an affine component of $\gb$.  By Remark \ref{breakpoints}, $(a,c)$
is therefore contained in an affine component of $\ga^{\gb}$, and
hence the slope of $\ga^{\gb}$ on $(a,c)$ is the leading slope of
$\ga^{\gb}$ on $A$, since the leading affine component of $\ga^{\gb}$
with non-trivial intersection with $A$ must now contain $(a,c)$.  The
slope of $\ga^{\gb}$ on $(a,c)$ is the product of the slope of
$\gb^{-1}$ on $(a,c)$, the slope of $\ga$ on $(a',c')$, and the slope
of $\gb$ on $(a',c'')$.  But the functions $\gb$ and $\gb^{-1}$ are
inverse, so if $x$ is in an affine component of $\gb^{-1}$ with slope
$1/s$, then $y=x\gb^{-1}$ is in an affine component of $\gb$ with slope
$s$.  Now then the image of $(a,c)$ under $\gb^{-1}$ must be in an
affine component $(a',d)$ of $\gb$, and the slopes of $\gb^{-1}$ on
$(a,c)$ and $\gb$ on $(a',d)$ are multiplicative inverses.  But now we
see that $(a',c'')\subset(a',d)$, so in particular, the lead slope of
$\ga^{\gb}$ on $(a,b)$ is the lead slope of $\ga$ on $(a',b')$.

The argument for the trailing slopes is similar. 
\qquad$\diamond$

Finally, we expand on the third point of Lemma \ref{finiteOrbitals}
to produce the following useful lemma.
\bl
\label{transitiveOrbital}
If $H\leq\ploi$ and $A=(a,b)$ is an orbital for $H$, then given any points
$c$, $d\in A$, with $c<d$, there is an element $g\in{}H$ so that
\mbox{$cg>d$}.
\el

pf: $[c,d]$ is contained in an orbital of $H$, and therefore it is
contained in the union of the orbitals of the elements of $H$.  Since
$[c,d]$ is compact, it follows that it is covered by a finite
subcollection $\mathscr{C}'$ of the orbitals of the elements of $H$.
This implies that there is a smallest positive integer $n$ and a
collection of signed orbitals $\mathscr{C} = \left\{(A_i,g_i)|1\leq
i\leq n, i\in\N\right\}$ of $H$ whose orbitals cover $[c,d]$.  If $n
=1$ we are done by Remark \ref{transitiveElementOrbital}, so assume
$n>1$.  Let us assume that the indexing of $\mathscr{C}$ respects the
left to right order of the orbitals in $\mathscr{C}$ (so that if $A_i$
and $A_j$ are two orbitals in $O_{\mathscr{C}}$, with $i<j$, then
given any point $y\in A_j$ there is $x\in A_i$ with $x<y$).  We note
in passing that by the minimality of $n$ the point $c$ is an element
of $A_1$ but not of $A_2$ and $d$ is an element of $A_n$ but $d$ is
not an element of $A_{n-1}$.

Improve $\mathscr{C}$ by supposing we chose signatures intelligently,
so that each $(A_i,g_i)\in\mathscr{C}$ satisfies the property that
$g_i$ moves points to the right on $A_i$.

For each $i \in \left\{1,2,\ldots,n\right\}$, let $B_i = A_i\cap
A_{i+1}$.  We note that for each such $i$, $B_i$ is a non-empty open
interval. Let $\delta_1$ be the minimum length of these intervals
$B_i$.  Now each orbital $A_i$ has length at least $\delta_1$.  Let
$\delta_2$ be the minimum of the two distances, one from the left end
of $A_1$ to $c$, and the other from $d$ to the right end of $A_n$.
Now let $\delta = \min(\delta_1,\delta_2)/2$.  It is immediate by
construction that each orbital of $\mathscr{C}$ has length at least
$2\delta$, and also that each interval $B_i$ created above has length
at least $2\delta$.  For each $k\in\left\{1,2,\ldots,n\right\}$, let
$x_k$ be a point in $A_k$ a distance less than $\delta$ from the left
side of $A_k$.  By Remark \ref{transitiveElementOrbital}, for each
integer $k$ where $1\leq k\leq n$, there is a positive integer $m_k$
so that $g_k^{m_k}$ will take $x_k$ to a point $y_k$ within $\delta$
of the right end of $A_k$.  Observing that for all integers $k$ where
$1\leq k<n$ we have that $y_k>x_{k+1}$, we see that $g =
g_1^{m_1}g_2^{m_2}\cdots g_n^{m_n}$ moves $x_1 $ to the right of
$y_n$.  But now $x_1<c<d<y_n$ by construction, so
$cg>d$.\qquad$\diamond$

Given a subgroup $G$ of $\ploi$ which has an orbital $A$, we will
often work with the image group of the homomorphism $\phi_A:G\to\ploi$
defined by the rule $g\mapsto g_A$, where $g_A$ is the unique element
of $\ploi$ which agrees with $g$ over the domain $A$, and behaves as
the identity elsewhere.  We will occasionally denote this image group
as $G_A$ (following Brin in \cite{BrinU}), and we will refer to it as
the \emph{projection of $G$ on the orbital $A$}.

\subsection{Thompson's group $F$ and balanced subgroups of $\ploi$}
Brin showed in \cite{BrinU} the following theorem:

\bt [Ubiquitous F]
\label{UbiquitousF}
If a group $H \leq \ploi{}$ has an orbital $A$ so that some element
$h\in{}H$ realizes one end of $A$, but not the other, then $H$ will
contain a subgroup isomorphic to Thompson's group $F$.
\et

The condition on the orbital is weak enough that one readily observes
that ``$\ploi{}$ is riddled with copies of $F$'', quoting Brin.  Hence
it becomes a natural question to ask what can be said about subgroups
of $\ploi{}$ which have no orbitals satisfying the ubiquity condition.

  We will say that an orbital $A$ of a group $H \leq \ploi{}$ is
\emph{imbalanced}\index{orbital!imbalanced} if some element $h \in H$
realizes one end of $A$, but not the other, and we will say $A$ is
\emph{balanced}\index{orbital!balanced} if whenever an element $h \in
H$ realizes one end of $A$, then $h$ also realizes the other end of
$A$ (note that $h$ might do this with two distinct orbitals).
Extrapolating, given a group $H \leq \ploi{}$, we will say that
\emph{$H$ is balanced}\index{group!balanced} if given any subgroup
$G\leq H$, and any orbital $A$ of $G$, every element of $G$ which
realizes one end of $A$ also realizes the other end of $A$.
Informally, $H$ has no subgroup $G$ which has an orbital that is
``heavy'' on one side.  In the case where $H$ has a subgroup $G$ with
an imbalanced orbital, then we will say that \emph{$H$ is
imbalanced}\index{group!imbalanced}.

\brk 

If $H\leq\ploi$ and $H$ is imbalanced, then $H$ has a subgroup
isomorphic to Thompson's group $F$.

\erk

Since $F'$ is non-trivial and simple (\cite{CFP}, Theorem 4.3), $F$ is
not solvable.  Thus imbalanced groups are not solvable.

The dynamics of balanced groups are much easier to understand than
those of imbalanced groups.  We will trade heavily on this in the
remainder.

\subsection{More on towers\label{towers}}
We now expand on the definitions of depth of a subgroup of $\ploi$
and height of any associated tower.

Suppose $G$ is a subgroup of $\ploi$, and $T$ is a tower of $G$.  If
there is an order preserving injection from $\N$ to $T$, then we will
say $T$ is \emph{tall}\index{tower!tall}, and if there is an order
preserving injection from the negative integers to $T$, then we will
say $T$ is \emph{deep}\index{tower!deep}.  If $T$ is both deep and
tall then we will say $T$ is a \emph{bi-infinite
tower}\index{tower!bi-infinite}, and we note that there will be an
order preserving injection from the integers to the tower.  We will
occassionally replace an infinite tower of one of these three types
with the image of the implied injection without comment, so that we
might refer to a tall, deep, or bi-infinite tower as countable, and
refer to the ``next'' element, etc., when this will not effect the
result of an argument.

Given a group $G\leq \ploi$, and an orbital $A$ of an element $g\in
G$.  We will define the \emph{depth of $A$ in
$G$}\index{orbital!depth} to be the supremum of the heights of the
finite towers which have their smallest element having the form
$(A,h)$ for some element $h\in G$.  If the depth of $A$ is an infinite
ordinal, we will say that $A$ is \emph{deep in
$G$}\index{orbital!deep}.  Symmetrically, we define the \emph{height
of $A$ in $G$}\index{orbital!height} to be the supremum of the heights
of the finite towers which have their largest element having the form
$(A,h)$ for some element $h\in G$.  If the height of $A$ is an
infinite ordinal, we will say that $A$ is \emph{high in
$G$}\index{orbital!high}.

The following are immediate from the definitions.
\brk
\label{trivialTower}

\be

\item Any subset of a tower is a tower.
\item If $T$ is a tower for some group $G\leq\ploi$, and $(A,g)\in T$,
and if $h$ is an element of $G$ with orbital $A$, then the set
$(T\backslash\left\{(A,g)\right\})\cup\left\{(A,h)\right\}$ is also a
tower for $G$.
\item If $T$ is a tower for some group $G\leq\ploi$, and $(A,g)\in T$,
and if $n$ is a non-zero integer, then
$(T\backslash\left\{(A,g)\right\})\cup\left\{(A,g^n)\right\}$ is also a
tower for $G$.
\item If $g\in G\leq\ploi$, and $T$ is a tower of $G$, and if $(A,g)$,
$(B,g)\in T$, then $A=B$.  That is, no signature appears twice in a tower.
\item Given a tower $T$ of a group $G\leq\ploi$, the group $H\leq G$
generated by the signatures of $T$ has an orbital $A$ which contains
all the orbitals of $T$.  
\item \label{conjugateTower} Given an element $k\in G\leq\ploi$ and any
tower $T$ for $G$, the set of signed orbitals
$T^k=\left\{(Ak,g^k)|(A,g)\in T\right\}$ is also a tower for $G$,
where the natural order of the signed orbitals of $T^k$ is equal to
the induced order from the signed orbitals of $T$.  
\ee 
\erk 

Given $k\in G\leq\ploi$ and a tower $T$ for $G$, the tower $T^k$
induced from the tower $T$ by the action of $k$ as discussed in item
\ref{conjugateTower} of Remark \ref{trivialTower} will be called the
\emph{tower conjugate to $T$ by the action of
$k$}\index{tower!conjugate}.  We will also say that the towers are
conjugate towers.  Conjugacy of towers for $G$ is an
equivalence relation on the set of towers for $G$.

Towers, as defined, are easy to find, but can be difficult to work
with.  For an arbitrary tower $T$, there are no guarantees about how
other orbitals of signatures of the elements of $T$ cooperate with the
orbitals of the tower.  We say a tower $T$ is an \emph{exemplary
tower}\index{tower!exemplary} if the following two additional
properties hold: \be

\item Whenever $(A,g)$, $(B,h)\in T$ then $(A,g)\leq(B,h)$ implies the
orbitals of $g$ are disjoint from both ends of the orbital $B$.

\item Whenever $(A,g)$, $(B,h)\in T$ then
$(A,g)\leq(B,h)$ implies no orbital of $g$ in $B$ shares an end with $B$.

\ee

The following lemmas indicate the plethora of
exemplary towers in $\ploi$, and will be used repeatedly.

\bl
\label{imbalancedTower}

Suppose $H\leq\ploi$, and that $G\leq H$ has imbalanced orbital $A =
(a,b)$.  Then $H$ admits an exemplary bi-infinite tower $E$ whose
orbitals are all in $A$.

\el

pf:

A short incorrect proof is that $H$ contains a copy of $F$ and $F$
contains an exemplary bi-infinite tower.  The problem is that Theorem
\ref{UbiquitousF} guarantees a copy of $F$, but does not guarantee
that its generators have one orbital each.

If we find an exemplary bi-infinite tower $E'$ for the projection
$G_A$ of $G$ on $A$, then by replacing the signatures of $E'$ with
elements of $G$ which agree with the signatures of $E'$ on $A$, we can
build a new exemplary tower $E$ for $G$ with the same orbitals.  Thus,
we may assume for the purposes of this argument that $G$ only has
orbital $A$.

  The remainder of the proof naturally breaks down into four stages:
\be
\item Building a bi-infinite tower near the right end of $A$.

Since $A$ is imbalanced for $G$, there is $g_0\in G$ so that $g_0$ has
an orbital $B_0$ which shares an end with $A$, but $g_0$ does not
realize the other end of $A$.  We will assume that $B_0$ shares its
right end with the right end of $A$, in particular, $B_0 = (a_0,b)$
for some $a_0\in A$, and there is $w\in A$ so that if $x\in \supp{g_0}$
then $w<x$.  We will further assume that $g_0$ moves points to the
right on the orbital $B_0$, so that all conjugates of $g_0$ move
points to the right on their corresponding orbitals.

By Lemma \ref{transitiveOrbital} there is an element $\alpha\in G$ so
that $r = w\alpha>a_0$.  Let $g_k = g_0^{(\alpha^{-k})}$ for all
integers $k$.  By construction, given any integer $k$, the element
$g_k$ has rightmost orbital of the form $(a_k,b)$, where
$a_{k+1}<\inf(\supp{g_k})$.

\item Thickening the tower.

We will construct a new sequence from the $g_k$ so that the righthand
orbitals of the new sequence do not get arbitrarily small, while
preserving all of the nesting properties we will need later.  Fix
$l\in (a_0, r)$.  

We now have $a<w<a_0<l<r<b$.  Now for each non-negative integer $k$,
we already have that $a_k<l$, so define $h_k = g_k$ when $k\geq 0$.
Since $h_0$ is already defined, we can define $h_{-k}$ inductively for
each positive integer $k$ as the first conjugate of $g_{-k}$ by
$h_{-k+1}^{-1}$ which has the property that the rightmost orbital
$(c_{-k},b)$ of $h_{-k}$ contains $(l,b)$.  For each integer $k$, let
$c_k$ represent the left end of the rightmost orbital of $h_k$.  In
particular, we have now defined a bi-infinite sequence of functions
$(h_k)_{k \in \Z}$ that satisfies the following properties:

\be

\item Given $k\in\Z$, the rightmost orbital of $h_k$ is $(c_k,b)$.

\item Given $k\in\Z$, $(l,b)\subset (c_k,b)$.

\item Given $k\in\Z$, $c_{k+1} < \inf(\supp{h_k})$.

\ee

\item Chopping off the right ends.

Note that the support of the element $g_{-1}$ lies to the right of
$r$.  For all $k\in\Z$, define $u_k = h_kg_{-1}^{-1}$.  Since the
$h_i$ are all conjugates of $g_0$, and $g_{-1}$ is a conjugate of
$g_0$, we see that the trailing slopes of all of the $h_k$ are the
same.  In particular, for each integer $k$, the fixed set of $u_k$ has
a component of the form $[e_k,b)$ where $l<r<e_k<b$.  Furthermore, for
each integer $k$, $u_k$ has an orbital $C_k$ of the form $(c_k,d_k)$
where $r<d_k\leq e_k<b$.  There may be other orbitals of $u_k$ to the
right of $d_k$ and the $d_i$ may not be ordered with respect to $i$.

\item Nesting the right ends (while keeping the left ends nested).

Let $v_0 = u_0$, and inductively define, for each positive integer
$k$, an element $v_k = u_k^{h_k^{m_k}}$ where $m_k$ is the smallest
positive integer so that the support of $v_{k-1}$ is fully contained
in the orbital $D_k = (c_k,r_k)$ of $v_k$ induced from
$C_k=(c_k,d_k)$.  Note that such an integer $m_k$ will always exist,
since $h_k$ moves points to the right on its rightmost orbital
$(c_k,b)$, and by reference to lemma \ref{transitiveElementOrbital}.

Now, for each positive integer $k$, inductively define elements
$v_{-k}$ by following the two step process below.

First, find the negative integer $n_k$ of smallest absolute value
which has the property that the closure of the support of $v_{-k}' =
u_{-k}^{(h_{-k}^{n_k})}$ is fully contained in the orbital $D_{-k+1}$
of $v_{-k+1}$ induced from $C_{-k+1}=(c_{-k+1},d_{-k+1})$.  At this
point $v_{-k-1}$ cannot be created since the orbital of $v_{-k}'$
induced from $C_{-k}$ might not contain the left end of the orbital
$(c_{-k-1},b)$ of $h_{-k-1}$.  

Now replace the elements of the sequences $(u_j)_{j\in\Z}$ and
$(h_j)_{j\in\Z}$ which have indices less than or equal to $-k$ by the
conjugate of each such element by $h_{-k}^{n_k}$.  Note $n_k<0$.  This
will do nothing to the element $h_{-k}$ of the sequence
$(h_j)_{j\in\Z}$, but all terms of $(h_j)_{j\in\Z}$ with $j<-k$ will
have their orbitals extended leftward by $|n_k|$ iterates of
$h_k^{-1}$.  Now define $v_{-k} = u_{-k}$.  At any inductive stage
$k$, note that by construction, all of the left ends of the orbitals
of the $h_j$, for indices $j<-k$, are inside the orbital $D_{-k}$, so
this inductive definition makes sense.  For each integer $k$, the
closure of the support of $v_{k-1}$ is a subset of this orbital $D_k$.
In particular, $E = \left\{(D_k,v_k)|k\in\Z\right\}$ is an exemplary
bi-infinite tower for $G$, and therefore for $H$, with all orbitals
$D_k$ in $A$, and where the index of $E$ respects the natural order on
the signed orbitals of $E$.  
\ee
\qquad$\diamond$

The following useful remark is easy to prove using the techniques
above, so we leave it to the reader.  This is a version of Lemma
\ref{imbalancedTChains}.

\brk

If $G$ is a subgroup of $\ploi$ that does not admit transition chains
of length two, then $G$ is balanced.

\erk

In the balanced group case, we have another way to sometimes find
bi-infinite, exemplary towers: 

\bl
\label{inconsistentTower}

Suppose $H$ is a balanced subgroup of $\ploi$ and $H$ has orbital $A$.
If there is $h\in H$ so that 

\be
\item $h$ moves points left on an orbital sharing an end of $A$, and 
\item $h$ moves points right on another orbital sharing the other end of $A$,
\ee 

then $H$ admits an exemplary bi-infinite tower $T$ with the orbitals
of \thinspace $T$ all contained in $A$.

\el
pf:

We will find such a tower for the projection of $H$ on $A$, then
replacing the signatures of our tower with elements of $H$ which agree
with our signatures over $A$ will create an exemplary tower for $H$
whose orbitals are all in $A$.  Thus, we shall assume that
$H$ has only the orbital $A$ for our discussion below.

First off, let us write $A=(a,b)$, so that we can refer to the ends
of $A$ in the argument below.

Let $h\in H$ so that $h$ satisfies the two points of the lemma
statement.  In particular, the two orbitals of $h$ specified by the
lemma statement are the first and last orbitals of $h$.  By
replacing $h$ with its inverse, if necessary, we can assume that $h$
moves points to the left on its first orbital and moves points to the
right on its last.  Let $F_h=\fix{h}\cap A$ represent the (non-empty)
fixed set of $h$ in $A$.

Let $r= \min(F_h)$ and $s = \max(F_h)$.  Note that since $F_h$ is
compact, both $r$ and $s$ exist.  Now by Lemma
\ref{transitiveOrbital}, there is $g\in H$ so that $rg>s$.  Let
$t$, $u\in A$ so that $t = rg$, and $u = sg$.  We now have
\[
a<r<s<t<u<b.
\]
Now Lemma \ref{transitiveElementOrbital} implies there is an integer
$k>0$ so that $th^k>u$.  Let us consider the element $g_0 =
h^{-k}gh^kg^{-1}$.  It is an immediate consequence of the second
point of Lemma \ref{breakpoints} that $g_0$ is the identity near the
ends of $A$.  A second valuable fact is that $[r,s]g_0\cap
[r,s]=\emptyset$.  To see this, let us compute the trajectory of a point
$x_0\in[r,s]$ under the action of $g_0$.  We will compute this
trajectory in stages using the defining product of $g_0$.  Each claim
to follow comes directly from our choice of $k$ and stated properties
of $g$.  Firstly, $x_0h^{-k}=x_1\in[r,s]$.  Secondly, $x_1g
=x_2\in [t,u]$.  Thirdly, $x_2h^k =x_3\in (u,b)$. Finally,
$x_3g^{-1} = x_4>s$.

Let $B_0$ represent the orbital of $g_0$ that contains $[r,s]$.  Note
by our construction of $g_0$ that the closure of the support of $g_0$
is contained in $A$ and that $g_0$ moves points to the right in $B_0$.
There is a smallest positive integer $n_1$ so that the orbital $B_1$
of $g_1 = g_0^{(h^{n_1})}$ induced from $B_0$ by the action of
$h^{n_1}$ contains the closure of the support of $g_0$, since repeated
conjugation of $g_0$ by $h$ increases the size of the orbital induced
from $B_0$ so that the ends of this induced orbital approach the ends
of $A$.  In particular, we can inductively define a bi-infinite
sequence $((B_i,g_i))_{i\in\Z}$ of signed orbitals of $H$ by the
property that $g_i = g_{i-1}^{(h^{n_i})}$ where $n_i$ is the smallest
positive integer so that the closure of the support of $g_{i-1}$ is
fully contained in the orbital $B_i$ induced from the orbital
$B_{i-1}$ by the action of $h^{n_i}$. (More formally, for negative
integers $i$, one defines $g_i$ from $g_{i+1}$ by saying that $g_i =
g_{i+1}^{\left(h^{-n_{i+1}}\right)}$ where $-n_{i+1}$ is the largest
negative integer so that the closure of the support of $g_i$ is
contained in $B_{i+1}$.  This makes sense since repeated conjugation
by $h^{-1}$ moves points in the first and last orbitals of $h$
arbitrarily close to the closed interval $[r,s]$).

By construction, the set $T= \left\{(B_i,g_i)|i\in\Z\right\}$ is an
exemplary bi-infinite tower for $H$ with all orbitals in $A$.
\qquad$\diamond$

In the remainder, given a subgroup of $\ploi$ with orbital $A$, and an
element $h$ so that the statement of the last lemma applies, we will
say that \emph{$h$ realizes the orbital $A$ inconsistently}.  If $h$
instead realizes both ends of $A$, but $h$ moves points in the same
direction near both ends of $A$ in $A$, then we will say that \emph{$h$
realizes $A$ consistently}.  Note that this does not imply that $h$ has
an empty fixed set in $A$.

We will use the previous two lemmas to show the following lemma, which
is a strengthening of Lemma \ref{tChainTowers}:

\bl 
\label{transitionTower}

Suppose $G$ is a subgroup of $\ploi$ and $G$ admits a transition
chain $\mathscr{C}$ of length two whose orbitals are contained in an
orbital $A$ of $G$, then $G$ admits an exemplary bi-infinite tower $T$
whose orbitals are all contained in $A$.

\el 

pf:

Suppose $\mathscr{C} = \left\{(A_1,g),(A_2,h)\right\}$, where $A_1 =
(a_1,b_1)$ and $A_2 = (a_2,b_2)$, and where $a_1<a_2<b_1<b_2$, and let
$O=(l,r)$ be the orbital of the group $H_1 = \langle g,h\rangle$ that
contains $(a_1,b_2)$.  Note in passing that $O\subset A$.  Now, all of
the points in $O$ are moved by at least one of $g$ and $h$, in
particular, we can assume without meaningful loss in generality that
$g$ realizes the end $l$ of $O$.  If $g$ does not also realize $r$
then $H_1$ is imbalanced, and we are done by application of Lemma
\ref{imbalancedTower}, therefore let us assume that $g$ also realizes
the end $r$.  If $g$ realizes the ends of $O$ inconsistently, then we
are done by an application of Lemma \ref{inconsistentTower}.
Therefore, let us assume that $g$ realizes the ends of $O$
consistently, and so we can further assume without meaningful loss of
generality that $g$ move points to the right on its first and last
orbitals in $O$.

Let $F_g$ represent the fixed set of $g$ in $O$.  The fixed set of $g$
in $O$ consists of finitely many closed, bounded intervals, all of
which are contained in the orbitals of $h$. In particular, by repeated
applications of Lemma \ref{transitiveElementOrbital}, we can find an integer
$n$ so that $h^n$ moves the fixed set of $g$ in $O$ completely off of
itself. (Note that the ends of the orbitals of $h$ are not in the
interiors of the components of the fixed set of $g$ in $O$.)  In
particular, the fixed set of $g$ is contained in the orbitals of the
element $k = [g,h^n]=g^{-1}h^{-n}gh^n$, while $k$ does not realize any
end of $O$.

We now will pass to consideration of the group $H_2 = \langle
g,k\rangle$, which has orbital $O$, with $g$ realizing both ends
consistently and having a nontrivial fixed set $F_g$ in $O$, and with
$k$ realizing neither end of $O$.  Let $x_1$, $y_1$ be elements of $O$ so
that if $z$ is any element of $F_g$ then $x_1<z<y_1$.  By Lemma
\ref{transitiveOrbital}, there is an element $q'$ of $H_2$ so that
$x_1q'>y_1$.  We can write $q'$ as a finite product of conjugates of $k$ by
powers of $g$, followed by $g^m$ for some integer power $m$.  Now the
element $q = q'g^{-m}$ satisfies $x_1q = v$ where $v$ is strictly to the
right of $F_g$, and also that $q$ is the identity near the ends of
$O$.  The first property just mentioned means that $F_g$ is contained
in a single orbital $C=(a_3,b_3)$ of $q$.

Let $x_2 = \min(\supp{q}\cap O)$ and $y_2 = \max(\supp{q}\cap O)$.  By
\ref{transitiveElementOrbital} there is a positive integer $i$ so that
$a_3<x_2g^i<b_3$ and so that $b_4 = b_3g^i>y_2$.  Let $p = q^{g^i}$.
The group $H_3 = \langle q, p \rangle$ has orbital $(a_3,b_4)$, and
the element $q$ realizes the end $a_3$ but not the end $b_4$, so that
this group is imbalanced, and therefore admits an exemplary
bi-infinite transition chain all of whose orbitals are contained in
$O$, and therefore in $A$.  But $H_3$ is a subgroup of our original
group $G$, so that $G$ admits this same bi-infinite exemplary
transition tower.  \qquad$\diamond$

In the following recall that $S_T$
is the set of signatures of a tower $T$.  

\bl

\label{TowersInCleanGroups}
If $G$ is a subgroup of $\ploi$ which does not admit transition chains of
length two and $G$ has a tower $T$, then $T$ is exemplary.

\el
pf: 

Suppose $G$ is a subgroup of $\ploi$ that does not admit
transition chains of length two, and $T$ is a tower for $G$.

Since $T$ is already a tower, we only need to show the following two
properties: 

\be
\item whenever $(A,g)$, $(B,h)\in T$ then $(A,g)\leq(B,h)$ implies the
orbitals of $g$ are disjoint from both ends of the orbital $B$, and

\item whenever $(A,g)$, $(B,h)\in T$ then $(A,g)\leq(B,h)$ implies no
orbital of $g$ in $B$ shares an end with $B$.  
\ee

The first property follows immediately since its violation would
immediately imply the existence of a transition chain of length two.

The lemma is obviously true if $T$ has height zero or one.  We now
assume that $T$ has at least two elements.  In particular, we can find
$(A,g)$ and $(B,h)\in T$ with $(A,g)\leq (B,h)$.  If the second
property is not satisfied, then we can chose our signed orbitals
intelligently so that we can find an orbital $A_1$ of $g$ in $B$ which
shares an end with $B$.  Let us suppose we have such an orbital $A_1$,
and derive a contradiction.

From the definition of the partial order on the signed orbitals of
$G$, and from the definition of a tower of $G$, we see that $A_1$
cannot extend all the way across $B$.  Therefore, and since our group
$G$ is balanced, we must have that some other orbital $A_2$ of $g$ in
$B$ shares the other end of $B$.  Let us assume without meaningful
loss of generality that $A_1$ is to the left of $A_2$.  In particular,
there are elements $a$, $a^*$, $b$, and $b^*\in [0,1]$ so that $B
= (a,b)$, $A_1 = (a,b^*)$, and $A_2 =(a^*,b)$.  By Lemma
\ref{transitiveElementOrbital}, there is an integer $n$ so that
$b^*h^n>a^*$.  Let $k= g^{h^n}$.  $k$ has an orbital $(a,b^*h^n)$,
while $g$ has an orbital $(a^*,b)$, where $a<a^*<b^*h^n<b$,
showing that $G$ admits transition chains of length two.  Therefore,
$T$ must also satisfy the second property above, and so $T$ is an
exemplary tower.  \qquad$\diamond$

\bc
\label{infiniteExemplaryTower}

If $G$ is a balanced subgroup of $\ploi$ and $G$ admits a tall tower in
some orbital $A$, or $G$ admits a deep tower in some orbital $A$, or $G$
admits a bi-infinite tower in some orbital $A$, then $G$ admits an
exemplary tall tower in $A$, or $G$ admits an exemplary deep tower in
$A$, or $G$ admits an exemplary bi-infinite tower in $A$,
respectively.

\ec 

pf: 

This follows from Lemma \ref{transitionTower}, and the details of the
proof of Lemma \ref{TowersInCleanGroups}, as below.  Suppose $T$ is an
infinite non-exemplary tower for a balanced subgroup $G$, where all
the orbitals of $T$ are contained in the orbital $A$ of $G$.  Since
$T$ is not exemplary, then we can produce a non-exemplary subtower $P
= \left\{(A_1,g_1), (A_2,g_2)\right\}$ of $T$ where $A_1\subset A_2$.
Suppose $B$ is the orbital of $G_P = \left<g_1,g_2\right>$ that
contains $A_1$.  Since $P$ is not exemplary, we must have that some
orbital of $g_1$ contains an end of $A_2$ or shares an end of $A_2$ in
$A_2$.  In the first case, $G_P$ admits a bi-infinite exemplary tower
in $B\subset A$ by Lemma \ref{transitionTower}, which has a subtower
of the appropriate type.  In the other case, $A_2 = B$, so that $A_2$
is actually an orbital of $G_P$.  Now since we are in the second case,
we must have that $g_1$ has some orbital $C$ in $A_2$ that shares one
end of $A_2$ in $A_2$.  Since $G_P$ is balanced, we have that $g_1$
realizes both ends of $A_2$ from within $A_2$.  Therefore, $g_1$ must
have another orbital $D$ in $A_2$ sharing the other end of $A_2$.  As
in the conclusion of the argument for Lemma \ref{TowersInCleanGroups},
some integer $k$ exists so that the conjugate $g_1^{g_2^k}$ will have an
induced orbital $Cg_2^k$ that will have non-trivial overlap with $D$,
so that we can form a transition chain of length two whose orbitals
are in $A_2$, and therefore in $A$.  Again by Lemma
\ref{transitionTower}, we can find an exemplary bi-infinite tower for
$G_P$ whose orbitals are in $A_2$, and this tower will have a subtower
of the appropriate type for $G$.  \qquad$\diamond$

The following remark is left to the reader.

\brk

If $n\in\N$, $G\leq\ploi$, and $G$ has depth $n$, then
\be
\item $G$ is balanced.
\item $G$ does not admit transition chains of length two.
\item $G$ only admits exemplary towers.
\ee  

\erk

\section{Geometric classification of solvable groups in $\ploi$}
We are now in a position to produce algebraic results by using our
geometric tools.  

The following two lemmas, and their proofs, explain why solvability
and depth are connected.  Together they complete the proof of Theorem
\ref{geoClassification} that a subgroup $H$ of $\ploi$ has derived
length $n$ if and only if the depth of $H$ is $n$.  Note that here and
in the remainder we use $[a,b] = a^{-1}b^{-1}ab$.

\bl
\label{shortTower}

Suppose $G$ is a subgroup of $\ploi$, $n$ is a positive integer, and
$A$ is an orbital of $G$.  If $G$ has a tower of height $n$ whose
orbitals are contained in $A$, then $G'$ has a tower of height $n-1$
whose orbitals are also contained in $A$.

\el 

pf: Suppose $T=\left\{(A_i,g_i)|1\leq i \leq n, i\in\N\right\}$ is a
tower of height $n$ all of whose orbitals are contained in the orbital
$A$ of $G$ and whose indexing respects the order of the elements of
$T$.  By Lemmas \ref{imbalancedTower} and \ref{transitionTower}, and
\ref{TowersInCleanGroups}, we can assume that $T$ is an exemplary
tower.  This follows since if $G$ admits a transition chain of length
two whose orbitals are in $A$, or if $G$ has a subgroup $H$ with
imbalanced orbital $B\subset A$, then we can find a bi-infinite
exemplary tower for $G$ whose orbitals are in $A$, and pass to an
exemplary sub-tower of height $n$, otherwise, $T$ must already be
exemplary.

We note that the lemma is trivially true if $n = 1$, so let us assume
$n\geq 2$.

For each integer $i$ with $2\leq i \leq n$ there is a smallest
positive integer $n_i$ so that the subset of the support of $g_{i-1}$
which is in $A_i$ is fully contained in a single fundamental domain of
$g_i^{n_i}$.  Define $h_1 = g_1$, and for each integer $i$
with $2\leq i \leq n$, define $h_i = g_i^{n_i}$.  Now define the set
$E = \left\{(A_i,h_i)|1\leq i\leq n, i\in\N\right\}$, which is an
exemplary tower.  $E$ has a nice property: if $i$ is an integer in
$1\leq i < n$, then the supports of $h_i$ and the supports of
$h_i^{h_{i+1}}$ are disjoint in $A_{i+1}$, so that $A_i$ is an orbital
of $[h_i,h_{i+1}]$.  For each integer $i$ in $1\leq i <n$ define $v_i
= [h_i,h_{i+1}]$.  Noting that $v_i\in G'$, we see that
$\left\{(A_i,v_i)|1\leq i \leq n-1, i\in\Z\right\}$ is a tower of
height $n-1$ for $G'$ whose orbitals are all in $A$.  \qquad$\diamond$

\bl
\label{solvableClassification}

If $G$ is a subgroup of $\ploi$ of depth $n$ for some positive
integer $n$, then $G'$ is a subgroup of $\ploi$ of depth $n-1$.

\el 

pf: 

Suppose $G$ is a subgroup of $\ploi$ with depth $n$.  $G$ must
be balanced by Lemma \ref{imbalancedTower}, and $G$ must have no
transition chains of length two by Lemma
\ref{transitionTower}.

Since $n>0$, $G$ is not the trivial group, and so $G$ has at least one
orbital.  Let $A = (a,b)$ be an orbital of $G$.  Note that if 
\[
\Upsilon =
\left\{B\,|\,B\textrm { is an orbital of some element of G}, B\subset
A\right\},
\]
then $A = \cup_{B\in\Upsilon}B$.  But $G$ has no transition chains of
length two or the depth of $G$ would be infinite, so $A$ can be
written as a union of a maximal chain of properly nested orbitals of
elements of $G$.  Taking these orbitals, paired with appropriate
signatures, we create a tower $T$ whose height is bounded above by
$n$.  Let the height of $T$ be $m$, and let $T = \left\{(A_i,g_i)\,|\,
1\leq i\leq m, i\in\Z\right\}$.  But now, by construction, $A$ is the
union of the orbitals $A_i$, all of which are contained in $A_m$, so
that $A_m$ must actually be $A$, and no element $g\in G$ has an
orbital $B$ properly containing $A_m$.  Thus, $(A_m,g_m)$ is a signed
orbital of depth one, and the orbitals of $G$ are precisely the
orbitals of depth one for $G$. (There are infinite-depth subgroups of
$\ploi$ which have orbitals that are not realized by any element of
the subgroup.)

Let $g$, $h\in G$, and consider the element $[g,h] \in G'$.  Let
$\Gamma$ be the set of all orbitals of $g$ and $h$.  Suppose $[g,h]$
has $t$ orbitals, where $t$ is some positive integer, and let
$\left\{B_i|1\leq i\leq t, i\in\N\right\}$ be the orbitals of $[g,h]$
in left to right order.  Both $g$ and $h$ fix $a$, so the slope of the
leftmost affine component of $[g,h]$ that intersects $A$ is the
product of the slopes of the leftmost affine components of $g^{-1}$,
$h^{-1}$, $g$, and $h$ which nontrivially intersect $A$, which product
is one.  In particular, $[g,h]$ cannot realize $A$, so no orbital of
$[g,h]$ is an orbital of $G$, and so no tower for $G'$ contains an
orbital of depth one for $G$, and thusly, all towers of $G'$ can have
height at most $n-1$.

By the last paragraph, we see that the depth of $G'$ is at most $n-1$.
By Lemma \ref{shortTower} $G'$ has a tower of depth $n-1$, so the
depth of $G'$ is actually $n-1$.
\qquad$\diamond$

Lemmas 2.8, 2.10, and 2.11 indicate some conditions which imply that a
group admits infinite towers (and would therefore be non-solvable).
The question arises as to whether there are groups that admit towers of
arbitrary height, but which do not admit any of the hypothesies of the
lemmas above.  Avoiding the easy infinitely generated case, we have
the following open question.

\bq
Does there exist a finitely-generated non-solvable subgroup of $\ploi$
which does not admit transition chains of length two?
\eq

\section{Standard restricted wreath products in $\ploi$}
The following is commonly used without comment in the remainder.  The
proof is left to the reader.

\brk

\label{nestedOrbitals}
Suppose $G\leq\ploi$ and $G$ does not admit transition
chains of length two.  If $(A,g)$ and $(B,h)$ are signed
orbitals of $G$ and $A\cap B\neq \emptyset$ then either
$\overline{A}\subset B$, $\overline{B}\subset A$, or $A=B$.
\erk

Below, we define a condition that may be satisfied by pairs of
elements of $\ploi$, and an algebraic construction that uses pairs of
elements of $\ploi$.  We point out that the construction produces
controlled results if the pair of elements involved satisfy the condition.
This last fact will be central to our proof of Theorem 1.2.

Suppose $h$, $k\in\ploi$ and they satisfy the properties found above
about how element orbitals intersect for elements in a group without
transition chains (whenever $A$ is an orbital of $h$ and $B$ is an
orbital of $k$, and $A\cap B\neq \emptyset$, then either $A = B$, or
$\bar{A}\subset B$, or $\bar{B}\subset A$).  We will say that $h$ and
$k$ satisfy the \emph{mutual efficiency condition} if given any
orbital of $C$ of $h$ that properly contains an orbital of $k$, then
the support of $k$ in $C$ is contained in a single fundamental domain
of $h$ in $C$, and the symmetric condition that whenever $D$ is an
orbital of $k$ that properly contains an orbital of $h$, then the
support of $h$ in $D$ is contained in a single fundamental domain of
$k$ in $D$.

Now let us consider the following algebraic commutator operation.

\newtheorem{construction}{Construction}

\begin{construction}[Double Commutator Operation]
Given two elements $h$, $k\in\ploi$ we can construct a third element
$[[h,k],k]$, which we will refer to as the \emph{double
commutator of $h$ and $k$}.

\end{construction}

Double commutators are nice to understand in the setting of a subgroup
of $\ploi$ with no transition chains of length two.  The following
useful lemma can be checked directly by the reader.

\bl
\label{dcFacts}

Let $h$, $k\in H$, where $H$ is a subgroup of $\ploi$ with no
transition chains of length two.  Suppose further that $h$ and $k$
satisfy the mutual efficiency condition.  If $f = [[h, k],k]$, then $f$
has the following properties:

\be

\item Every orbital of $h$ properly contained in an orbital of $k$ is
an orbital of $f$.

\item Every orbital of $f$ is properly contained in an orbital of $k$ that
contains (perhaps not properly) an orbital of $h$.

\ee 

\el

We are now ready to approach a proof of Theorem 1.2 classifying the
top group of any standard restricted wreath product embedded in
$\ploi$.  First, we need the following lemma.

\bl
\label{solveWreathStructure}

The group $\Z\wr(\Z\times\Z)$ fails to
embed in $\ploi$.

\el

pf:

Suppose that $\phi:\Z\wr(\Z\times\Z)\to \ploi$ is an embedding.  Let
$\hat{B}$ and $\hat{T}$ represent the base and top subgroups of
$\Z\wr(\Z\times\Z)$, respectively.  Now let $H =
\textrm{Image}(\phi)$, and let $B=\hat{B}\phi$ and $T = \hat{T}\phi$
represent the images of the base and top subgroups.  We will refer to
$B$ and $T$ as the base and top groups of $H$.  Note that since $H$ is
solvable, $H$ admits no transition chains of length two.

We first observe that if $b\in B$ and $t\in T$ are both non-trivial,
then $[[b,t],t]\in B$ and $[[b,t],t]\neq 1$.

The first point follows inductively from the fact that a commutator of
an element of the base group with an element of the top group in a
wreath product is an element of the base group.

The second point follows from the first point because of a property of
standard restricted wreath products; Non-trivial elements of the base
group never commute with non-trivial elements of the top group when
the top group is torsion free.  In our case $[[b,t],t]$ is a
commutator of a nontrivial element of the base group ($[b,t]\in B$ by
the previous point) with a non-trivial element of the top group ($t\in
T$), and hence is non-trivial.

Consider the elements $a = ((0)_{i\in\Z\times\Z}, (1,0))$ and $b =
((0)_{i\in\Z\times\Z}, (0,1))$ in $\Z\wr(\Z\times\Z)$.  Note that
$\hat{T} = <a,b>$.  Let $\alpha = a\phi$, and $\beta = b\phi$.  We now
have that $T=\langle \ga,\,\gb\rangle$, and any element of $T$ can be
written in the form $\ga^k\gb^j$ for some integers $k$ and $j$.

We will now find a non-trivial element $\gamma\in B$ which must
commute with $\alpha$, contradicting the nature of $H$ as a standard
restricted wreath product of non-trivial groups with top group torsion
free.  Our approach is to pick a candidate for $\gamma$ and then repeatedly
improve $\gamma$ until it will provide us with the desired contradiction.

First let $\gamma$ be any non-trivial element of $B$.  There are
integers $k$ and $j$ so that $\gb^k$ and $\gamma^j$ satisfy the mutual
efficiency condition. Replace $\gamma$ with the double commutator
$[[\gamma^j,\gb^k],\gb^k]$.  Note that the new $\gamma$ is a
non-trivial element of $B$ and that the closure of the support of
$\gamma$ now lies entirely in the orbitals of $\gb$.

If $\ga$ and $\gb$ do not have disjoint supports, then any connected
component of the intersection of their supports is actually an orbital
of each element (or $\ga$ and $\gb$ will fail to commute).  Brin and
Squier show in \cite{picric} that if two one-orbital functions in
$\ploi$ share their orbital and commute then the group they generate
has an element which is a common root of these two generators.  As a
consequence, if $A$ is an orbital shared by both $\ga$ and $\gb$ then
there is an element in $T$ which behaves as a root of both
$\ga$ and $\gb$ over $A$.  In particular, there are non-zero integers
$m$ and $n$ so that the non-trivial product $\theta = \ga^m\gb^n$
behaves as the identity over $A$.  We can now find new integers $j$
and $k$ so that $\gamma^j$ and $\theta^k$ satisfy the mutual
efficiency condition.  Let us replace $\gamma$ by the element
$[[\gamma^j,\theta^k],\theta^k]$, which is now a non-trivial element
of the base group which can only have orbitals in the same orbitals of
$\gb$ as before, but which will have no orbitals in $A$.

Repeating the process above for each of the orbitals shared by $\ga$
and $\gb$ will produce a non-trivial element $\gamma$ of the base
group whose support is disjoint from the support of $\ga$.  This
element must commute with $\ga$.  But non-trivial elements of the base
group of a standard restricted wreath product cannot commute with
non-trivial elements of the top group when the top group is torsion
free, so $\Z\wr(\Z\times\Z)$ cannot embed in $\ploi$.

\qquad$\diamond$

Now we are ready to prove Theorem 1.2, namely, any standard restricted
wreath product of non-trivial groups embedded in $\ploi$ must have the
top group isomorphic with the integers.

\emph{Proof of Theorem \ref{wreathTopClassification}}:

Suppose $H\leq\ploi$ and $H\cong C\wr T$ where $C$ and $T$ are
non-trivial groups.

By Theorem 10.10.1 in P. H. Neumann's \cite{NeumannW}, any group which
can be written as a standard restricted wreath product of non-trivial
groups will have the top group determined up to isomorphism.

By the proof of Theorem 4.8 in \cite{CFP}, any non-abelian subgroup of
$\ploi$ must contain a copy of $\bigoplus_{i\in\Z}\Z$.  In particular,
if the top group of $H$ is non-abelian, then it contains a copy of
$\Z\times\Z$.  By taking a specific non-trivial element of the bottom
group, and passing to the subgroup of the bottom group it generates (a
group isomorphic with $\Z$) we can find a subgroup of $H$ which is
isomorphic with $\Z\wr(\Z\times\Z)$.  This last group does not embed
in $\ploi$, so we must have that the top group is abelian.

It is a consequence of theorem 5.4 in Brin and Squier's paper
\cite{picric} that no element of $\ploi$ has infinitely many
roots.  Therefore, the top group is a torsion free infinite abelian
group which does not admit $\Z\times\Z$ as a subgroup, and for which
no element has infinitely many roots.  In particular, the top group
must be isomorphic with the integers.  \qquad$\diamond$

\setlength{\baselineskip}{.67\baselineskip}
\bibliographystyle{amsplain}
\bibliography{ploiBib}
\end{document}